\documentclass[amsthm,thmtools]{ifacconf}

\usepackage{natbib}
\usepackage{amsmath,amssymb,amsfonts,graphicx,url}

\usepackage{times}
\usepackage{dsfont}
\usepackage{mathtools}
\usepackage[utf8]{inputenc}
\usepackage[T1]{fontenc}

\usepackage{graphicx,color,fancyhdr}
\usepackage{enumerate}
\usepackage{multirow}

\newcommand{\Expect}{\mathbb{E}}
\newcommand{\ud}{\,\mathrm{d}}

\renewcommand{\tilde}{\widetilde}

\newcommand{\real}{\mathbb{R}}

\newcommand{\calZ}{\mathcal{Z}}
\newcommand{\sfK}{\mathsf{K}}
\newcommand{\Ricc}{\text{Ricc}}
\newcommand{\Xbar}{\bar{X}}
\newcommand{\Bbar}{\bar{B}}
\newcommand{\mbar}{\bar{m}}
\newcommand{\Sigmabar}{\bar{\Sigma}}
\newcommand{\Wbar}{\bar{W}}

\newcommand{\mN}{m^{(N)}}
\newcommand{\SigN}{\Sigma^{(N)}}
\newcommand{\calU}{\mathcal{U}}
\newcommand{\calV}{\mathcal{V}}
\newcommand{\calR}{\mathcal{R}}
\newcommand{\calQ}{\mathcal{Q}}
\newcommand{\pibar}{\bar{\pi}}
\newcommand{\calW}{\mathcal{W}}

\newcommand{\newP}[1]{\medskip\noindent{\bf #1}}

\newtheorem{proposition}[thm]{Proposition}

\newtheorem{remark}{Remark}

\begin{document}

\begin{frontmatter}
	
	\title{Optimality vs Stability Trade-off in \\Ensemble Kalman Filters} 
	
	
	\author[AHT]{Amirhossein Taghvaei} 
	\author[PGM]{Prashant G. Mehta} 
	\author[TTG]{Tryphon T. Georgiou}
	
	\address[AHT]{Department of Aeronautics and Astronautics, University of Washington, Seattle 
		CA (e-mail: amirtag@uw.edu).}
	\address[PGM]{Department of Mechanical Science and Engineering, University of Illinois at Urbana-Champaign, 
		IL  (e-mail: mehtapg@illinois.edu)}
	\address[TTG]{Department of Mechanical and Aerospace Engineering, University of California, Irvine, 
	CA  (e-mail: tryphon@uci.edu)}

\begin{abstract}
	This paper is concerned with optimality and stability analysis of a family of ensemble Kalman filter (EnKF) algorithms. EnKF is commonly used as an alternative to the Kalman filter for high-dimensional problems, where storing the covariance matrix is computationally expensive. The algorithm consists of an ensemble of interacting particles driven by a feedback control law. 
	The control law is designed such that, in the linear Gaussian setting and asymptotic limit of infinitely many particles, the mean and covariance of the particles follow the exact mean and covariance of the Kalman filter. The problem of finding a control law that is exact does not have a unique solution, reminiscent of the problem of finding a transport map between two distributions. A unique control law can be identified by introducing  control cost functions, that are motivated by the optimal transportation problem or Schr\"odinger bridge problem. The objective of this paper is to study the relationship between optimality and long-term stability of a family of exact control laws. Remarkably, the control law that is optimal in the optimal transportation sense leads to an EnKF algorithm that is not stable.
\end{abstract}
	
	
\end{frontmatter}



\section{Introduction}
Consider the linear system 
\begin{subequations}\label{eq:model}
	\begin{align}
	\ud X_t &= AX_t + \sigma_B \ud B_t,\quad X_0 \sim N(m_0,\Sigma_0)\\
	\ud Z_t &= HX_t \ud t +  \ud W_t, 
	\end{align}
\end{subequations}
	where $X_t\in \real^d$ is the state of the system at time $t$, $Z_t \in \real^m$ is the observation process, $B_t \in \real^d$ and $W_t \in \real^m$ are mutually independent standard Brownian motions,   and $A$, $H$, $\sigma_B$ are matrices of appropriate dimensions. The distribution of the initial state $X_0$ is Gaussian $N(m_0,\Sigma_0)$ with mean $m_0$ and covariance $\Sigma_0$. 
	
	The filtering problem is to compute the {\em posterior distribution} of $X_t$ conditioned on the filtration generated by the history of observations $\calZ_t:=\sigma(Z_s;~0\leq s\leq t)$. 	
	For the linear system~\eqref{eq:model}, the filtering problem admits an explicit solution: the posterior distribution is Gaussian $N(m_t,\Sigma_t)$ with mean and covariance  governed by the Kalman-Bucy filter equations~\citep{kalman-bucy}: 
	\begin{subequations}\label{eq:KF}
		\begin{align}
		\ud m_t &= Am_t \ud t + \sfK_t(\ud Z_t - Hm_t\ud t )=:\mathcal T_t(m_t,\Sigma_t),\label{eq:KF-mean}\\
		\frac{\ud \Sigma_t}{\ud t} &= A\Sigma + \Sigma A^\top + \Sigma_B -  \Sigma H^\top H \Sigma  =:\Ricc(\Sigma_t),\label{eq:KF-var}
		\end{align}
	\end{subequations}
where $\sfK_t:=\Sigma_tH^\top$ is the Kalman gain and  $\Sigma_B:=\sigma_B\sigma_B^\top$.  The notation $\mathcal T_t(\cdot,\cdot)$ and $\Ricc()$ is used to  identify the update law for the mean and covariance respectively. 


Ensemble Kalman filter (EnKF) is a Monte-Carlo-based numerical algorithm that is designed to approximate the solution to the filtering problem~\citep{evensen1994sequential,whitaker2002ensemble,reich11,Reich-ensemble}. 
EnKF is widely used in in applications (such as weather prediction) where the state
dimension $d$ is very high; cf.,~\citep{Reich-ensemble,houtekamer01}.  The high-dimension of the
state-space provides a significant computational challenge {\em even}
in linear Gaussian settings.  For such problems, an EnKF
implementation may require less computational resources (memory and
FLOPS) than a Kalman filter~\citep{houtekamer01,evensen2006}.  

The EnKF algorithm is designed in two steps:
\begin{enumerate}
	\item[(i)] a controlled stochastic process $\Xbar_t$ is constructed, whose conditional distribution given $\calZ_t$  is equal to the posterior distribution of $X_t$:
	\begin{equation}\label{eq:exactness}
	\text{(exactness)}\quad P(X_t\in \cdot |\calZ_t) = 
	P(\Xbar_t\in \cdot |\calZ_t).
	\end{equation} 
	\item[(ii)] an ensemble of $N$ stochastic processes $\{X^i_t\}_{i=1}^N$ is simulated to empirically approximate the distribution of   $\Xbar_t$:
	\begin{equation}\label{eq:approx}
P(\Xbar_t\in \cdot |\calZ_t)\approx
	\frac{1}{N} \sum_{i=1}^N \delta_{X^i_t}(\cdot ).
	\end{equation} 
\end{enumerate}
The process $\Xbar_t$ is referred to as {\em mean-field process} and the stochastic processes $\{X^i_t\}_{i=1}^N$ are referred to as {\em particles}.  The property~\eqref{eq:exactness} is referred to as {\em exactness}.

\renewcommand{\arraystretch}{1.2}
\begin{table*}[t]\label{tab:Grqs}
	\caption{Description of  three established forms of EnKF.  $G_t,r_t,q_t$ are parameters of the mean-field process update law~\eqref{eq:Xbar}. Stability rate and approximation error are described in Prop.~\ref{prop:stability} and~\ref{prop:XN} respectively. }
	\centering
	\begin{tabular}{|c|c|c|c|c|c|c|} \hline 
		Algorithm &  Acronym & $G_t$ & $r_t$ & $q_t$ & Stability  rate & steady-state error \\\hline
		EnKF with perturbed observation~\citep[Eq. (26)]{reich11} &  P-EnKF &$A - \Sigmabar_tH^2$ & $\sigma_B$ & $\Sigmabar_tH$ & $\lambda_0$  & $\propto N^{-1}(\Sigma_B+ H^2)$\\
		Square-root EnKF~\citep[Eq. (3.3)]{Reich-ensemble} & S-EnKF & $A - \frac{1}{2}\Sigmabar_tH^2$ & $\sigma_B$ & $0$ & $\frac{1}{2}(\lambda_0-A)$  & $\propto N^{-1}\Sigma_B$ \\ 
		Deterministic EnKF~\citep{AmirACC2016} & D-EnKF & $A - \Sigmabar_tH^2 + \frac{1}{2}\Sigma_B\Sigmabar^{-1}_t$ & $0$ & $0$ & $0$ & $	0$ \\\hline 
	\end{tabular}
\end{table*}
\renewcommand{\arraystretch}{1.0}

The motivation to study the EnKF algorithm is two-fold:

	(i) As mentioned above, EnKF algorithm is computationally efficient compared to Kalman filter algorithm for high-dimensional problems~\citep{Reich-ensemble}. The computational cost of EnKF scales as $O(Nd)$, wheres it scales as $O(d^2)$ for Kalman filter. 

	(ii) The two-step design procedure can be generalized to {\em nonlinear} and {\em non-Gaussian} setting. The result is the feedback particle filter (FPF) algorithm~\citep{taoyang_TAC12,yang2016}. Therefore, EnKF can be considered as special case of  FPF, and understanding EnKF is insightful for understanding the FPF algorithm. 

Here is the outline  and summary of the paper:

 (i) The problem of constructing a mean-field process $\Xbar_t$, such that exactness property~\eqref{eq:exactness} is satisfied, is addressed in Section~\ref{sec:mean-field}.  It is shown that exact mean-field process is not unique. The family of exact processes is identified in Prop.~\ref{prop:Xbart}. 
Three members of the family, that correspond to three established forms of EnKF algorithm, are presented in Table~\ref{tab:Grqs}. 

(ii) The stability of the mean-field process is studied in Section~\ref{sec:stability}.  It is shown that, the mean-field processes exhibit different stability behaviour.  In particular, the deterministic EnKF, that is constructed from limit of incremental optimal transportation maps, is not stable. While, other forms EnKF that involve stochastic terms, are stable. 

(iii) Finally, the  particle system and the analysis of the approximation error~\eqref{eq:approx} is presented in Section~\ref{sec:particles}. It is observed that the EnKF algorithms that were stable, exhibit larger approximation error, because of the presence of stochastic terms.   


This paper builds on the growing literature on the design and the analysis of the EnKF algorithms.  Three forms of the EnKF algorithm are of importance: EnKF with perturbed observation~\citep{evensen1994sequential,reich11}, square-root form of EnKF~\citep{whitaker2002ensemble,Reich-ensemble} which is the same as FPF algorithm with constant gain approximation~\citep{yang2016}, and deterministic EnKF~\citep{AmirACC2016,taghvaei2020optimal}. This paper is also related to the stability and tthe error analysis of the EnKF algorithm in the  discrete-time setting~\citep{gland2009,mandel2011convergence,tong2016nonlinear,stuart2014stability,kwiatkowski2015convergence}, and the continuous-time setting~\citep{delmoral2016stability,jana2016stability,delmoral2017stability,bishop2018stability}.



\newP{Assumption I:} It is assumed that the processes $X_t$ and $Z_t$ are scalar, i.e. $d=m=1$.  This assumption is made to simplify the exposition. The possible extension to vector-valued case is briefly discussed as a remark in each Section.

\section{Construction of mean-field process}\label{sec:mean-field}
In order to construct an exact mean-field process $\Xbar_t$, consider the following sde
\begin{align}
\nonumber 
\ud \Xbar_t = ~&\calU_t(\Xbar_t)\ud t + \calV_t(\Xbar_t) \ud Z_t +\calR_t(X_t) \ud \Bbar_t \\
\label{eq:Xbar-general}
+ &\calQ_t(\Xbar_t) \ud \Wbar_t,  \quad 
\Xbar_0 \sim \pibar_0,
\end{align}
driven by control laws $\calU_t(\cdot)$, $\calV_t(\cdot), \calR_t(\cdot), \calQ_t(\cdot)$, where $\Bbar_t$ and $\Wbar_t$ are independent copies of $B_t$ and $W_t$ of~\eqref{eq:model}, and $\bar{\pi}_0$ is the initial distribution. The following result characterizes control laws that lead to a mean-field process $\Xbar_t$ with exactness property~\eqref{eq:exactness}.
\begin{proposition}\label{prop:Xbart}
The mean-field process~\eqref{eq:Xbar-general} satisfies the exactness property~\eqref{eq:exactness} if  the initial distribution $\pibar_0$ is Gaussian $N(m_0,\Sigma_0)$, 
and 
\begin{subequations}\label{eq:UV}
\begin{align}
\calU_t(x) &= (A-\Sigmabar_tH^2)\mbar_t +   G_t(x-\mbar_t),\\
\calV_t(x) &= \Sigmabar_t H ,\\
\calR_t(x) &= r_t ,\quad
\calQ_t(x)= q_t, 
\end{align} 
\end{subequations}
where $\mbar_t=\Expect[\Xbar_t|\calZ_t]$, $\Sigmabar_t=\Expect[(\Xbar_t-\mbar_t)^2|\calZ_t]$, and $G_t,r_t,q_t\in \real$ satisfy 
\begin{equation}\label{eq:Grq-1}
2G_t\Sigmabar_t + r_t^2  + q_t^2 = \Ricc(\Sigmabar_t).
\end{equation}
\end{proposition}
Using the form of control laws in~\eqref{eq:UV}, the sde~\eqref{eq:Xbar-general} is 
\begin{align}\nonumber
\ud \Xbar_t =  ~&\mathcal T_t(\mbar_t,\Sigmabar_t) + G_t(X^i_t-\mbar_t)\ud t + r_t \ud \Bbar_t \\\label{eq:Xbar}
+ &q_t \ud \Wbar_t ,\quad \Xbar_0 \sim \pibar_0
\end{align}
where $\mathcal T_t(m,\Sigma)$ is defined in~\eqref{eq:KF}. The fact that $\Xbar_t$ satisfies the exactness property~\eqref{eq:exactness} is observed by noting that: (i) $\mbar_t=m_t$ and $\Sigmabar_t=\Sigmabar_t$, because the time-evolution of $\mbar_t$ and $\Sigmabar_t$ is identical to the Kalman filter equations~\eqref{eq:KF}; (ii) the distribution of $\Xbar_t$ is Gaussian because the sde~\eqref{eq:Xbar} is linear upon replacing  $\mbar_t=m_t$ and $\Sigmabar_t=\Sigmabar_t$ and the initial distribution $\pi_0$ is Gaussian. 

There are three established choices for the parameters $G_t,r_t,q_t$ that are tabulated in Table~\ref{tab:Grqs}, where each row corresponds to a specific  form of EnKF algorithm. 

\begin{remark}
	The form of the sde~\eqref{eq:Xbar-general} may not  be general enough to capture all possible stochastic processes~$\Xbar_t$ that achieve the exactness property~\eqref{eq:exactness}. The particular form of sde~\eqref{eq:Xbar-general} is motivated by its appealing control theoretic form, where all the control terms are assumed to be of feedback form. For a more general prescription of exact mean-field processes, see~\cite{abedi2019gauge}. 
\end{remark}

\begin{remark}
	For the vector-valued case, $\calU(x)$ involves additional divergence free term (e.g. $\psi(x)=\Omega\Sigmabar^{-1}(x-\mbar_t)$ with skew-symmetric matrix $\Omega$) that does not effect the distribution~\citep{taghvaei2020optimal}.  
\end{remark}

\begin{remark}
	The control law governing the deterministic EnKF in Table~\ref{tab:Grqs} is optimal with respect to control cost associated with optimal transportation problem. In particular,  it is obtained as the continuous-time limit of infinitesimal optimal transportation maps between the Gaussian distributions that are given by the Kalman filter~\citep{taghvaei2020optimal}. Also, the control law  governing the square-root EnKF in Table~\ref{tab:Grqs} is optimal with respect to a control cost associated with the Schr\"odinger bridge problem with prior dynamics given by~\eqref{eq:model}~\citep{chen2016Linear}.  
\end{remark}

\section{Stability of the mean-field process}\label{sec:stability}
Let $\pi_t$ denote the exact posterior distribution given by Kalman-Bucy filter~\eqref{eq:KF}, and $\bar{\pi}_t$ denote the distribution of the mean-field process given by~\eqref{eq:Xbar}. Proposition~\ref{prop:Xbart} informs that the exactness condition $\bar{\pi}_t=\pi_t$ is satisfied if the initial distribution $\bar{\pi}_0 = \pi_0$ and \eqref{eq:Grq-1} holds. 
 The objective is to study the error between $\bar{\pi}_t$ and $\pi_t$ if the initial distributions $\pibar_0$ and $\pi_0$ are not equal.


The convergence analysis is carried out by analyzing the convergence of $\pibar_t$ to the Gaussian distribution $\tilde{\pi}_t = N(\bar{m}_t,\bar{\Sigma}_t)$ with the same mean and variance as the mean-field process, and the convergence of the mean and variance to the mean and variance given by the Kalman-Bucy filter equations. 
The convergence result is presented in the following proposition. We use  $2$-Wasserstein metric~\citep{villani2003}, denoted by $\mathcal W_2(\cdot,\cdot)$, in order  to measure the error between distributions, and we make  the following assumption.

\newP{Assumption II:} The linear system~\eqref{eq:model} is controllable  and observable. In the scalar case, this amounts to $\sigma_B\neq 0$ and $H\neq 0$.

%
%
%


\medskip

\begin{proposition}\label{prop:stability}
	Consider the mean-field process \eqref{eq:Xbar} under assumption (I)-(II). 
Then, \\
	(i) The mean and the variance  of the mean-field process converge to the exact mean and variance given by the Kalman-Bucy filter~\eqref{eq:KF}:
				\begin{align}
						\Expect[|\mbar_t-m_t|^2] &\leq \text{(const.)}e^{-2\lambda_0 t}(|\mbar_0-m_0|^2 + |\Sigmabar_0-\Sigma_0|^2 ), \nonumber\\ 
						|\Sigmabar_t-\Sigma_t|&\leq \text{(const.)}e^{-2\lambda_0 t} |\Sigmabar_0-\Sigma_0|, 
					\end{align}
			where $\lambda_0 = (A^2 + H^2\Sigma_B)^{\frac{1}{2}}$.\\
	(ii) The error between the mean-field distribution $\bar{\pi}_t$ and the Gaussian distribution $\tilde{\pi}_t = N(\bar{m}_t,\Sigmabar_t)$ is bounded by:
			\begin{align}
		\mathcal W_2(\bar{\pi}_t,\tilde{\pi}_t) \leq e^{\int_0^t G_s \ud s}~{\calW}_2(\pibar_0,\tilde{\pi}_0)
	\end{align}
\\
	(iii) Combining part (i) and (ii), the total error between the mean-field distribution and the exact filter is bounded by 
			\begin{align}\label{eq:W2}
				\Expect[\mathcal W_2(\pi_t,\bar{\pi}_t)] \leq Ce^{-\lambda_0 t}\mathcal W_2(\tilde{\pi}_0,\pi_0) + e^{\int_0^t G_s \ud s}~{\calW}_2(\pibar_0,\tilde{\pi}_0)
			\end{align}
			where $C$ is a constant independent of time $t$.
%
%
%
%

%

\end{proposition}
The convergence result~\eqref{eq:W2} decomposes the error into two terms. 
The first term in the error is due to the incorrect specification of the initial mean and variance. In particular, \[\mathcal W_2(\tilde{\pi}_0,\pi_0)^2 = (\bar{m}_0-m_0)^2 + (\sqrt{\Sigmabar_0} - \sqrt{\Sigma_0})^2.\] The contribution from this term converges to zero with an exponential rate, independent of the choice made for $G_t,q_t,r_t$ in~\eqref{eq:Xbar} as long as~\eqref{eq:Grq-1} holds. The bound follows from the the stability of the Kalman-Bucy filter which holds under controllability  and observability of the linear system~\citep{ocone1996}. 
 
The second term in the error is due to the fact that the initial distribution is not Gaussian. It is controlled by the stability of the mean-field process~\eqref{eq:Xbar} which, in contrast to the first term in the error, depends on the choice for $G_t$. 
For three choices of $G_t$, determined by the three forms of the EnKF tabulated in Table~\ref{tab:Grqs}, the following holds 
\begin{subequations}\label{eq:W2-EnKF}
	\begin{align}
		\text{(P-EnKF)}\quad\quad \quad e^{\int_0^t G_s \ud s}&\leq   \text{(const.)}e^{-\lambda_0 t}, \\
		\text{(S-EnKF)}\quad \quad \quad  e^{\int_0^t G_s \ud s}   &\leq \text{(const.)}  e^{-\lambda_1 t}, \\
		\text{(D-EnKF)}\quad  \quad \quad e^{\int_0^t G_s \ud s} &= \sqrt{\frac{\bar{\Sigma}_t}{\bar{\Sigma}_0}},
	\end{align}
\end{subequations}
where $\lambda_1=\frac{\lambda_0-A}{2}$. Following conclusions are in order: 
\begin{enumerate}
	\item[(i)]  Both P-EnKF and S-EnKF  are stable, i.e. the error converges to zero as $t\to \infty$. Moreover, $\lambda_0>\lambda_1>0$ (because $\lambda_0 =   (A^2 + H^2\sigma_B^2)^{\frac{1}{2}} > |A|$). Therefore, the convergence rate of P-EnKF is strictly larger than  the convergence rate of S-EnKF.
	\item[(ii)] D-EnKF is not stable. If the initial distribution is non-Gaussian, it remains non-Gaussian. 
In fact, one can establish the asymptotic lower-bound
	\begin{equation*}
		\liminf_{t \to \infty} \mathcal W_2(\pibar_t,\pi_t) \geq \text{(const.)} {\calW}_2(\pibar_0,\tilde{\pi}_0) 
	\end{equation*}
	implying that the error remains positive if the initial distribution is not Gaussian.
\end{enumerate}

\begin{remark} The result~\eqref{eq:W2} can be extended to vector-valued case by replacing $\exp(\int_0^tG_s\ud s)$ with the state transition matrix $\Phi_t$ defined as the solution to $\frac{\ud}{\ud t}\Phi_t = G_t \Phi_t$ with $\Phi_0=I$. 
\end{remark}


\section{Particle system and error analysis}\label{sec:particles}
The system of particles $\{X^i_t\}_{i=1}^N$ is constructed from the mean-field process~\eqref{eq:Xbar} by empirically approximating the mean and covariance:
\begin{align}\nonumber
\ud X^i_t =  &\mathcal T_t(\mN_t,\SigN_t) + G_t(X^i_t-m^{(N)}_t)\ud t + r_t \ud B^i_t \\\label{eq:particles}
+ &q_t \ud W^i_t ,\quad X^i_0 = x^i_0,\quad \text{for}\quad i=1,\ldots,N
\end{align}
 where $\mN_t=N^{-1}\sum_{i=1}^N X^i_t$ and $\SigN_t=(N-1)^{-1}\sum_{i=1}^N(X^i_t-\mN_t)^2$ are the empirical mean and the empirical covariance of the particles respectively, and $\{B^i_t\}_{i=1}^N$ and $\{B^i_t\}_{i=1}^N$ are independent copies of $\Bbar_t$ and $\Wbar_t$.  The initial state of the particles are denoted by $\{x^i_0\}_{i=1}^N$. 
 
 The objective is to analyze the error between the mean-field distribution $\bar{\pi}_t$ and the empirical distribution of the particles 
\begin{equation*}
\pi^{(N)} _t(\cdot)= \frac{1}{N}\sum_{i=1}^N \delta_{X^i_t}(\cdot),
\end{equation*}
where $\delta_x$ is the Dirac delta distribution located at $x$. The result is presented for the convergence of the empirical covariance to the mean-field covariance under the following assumption:

\newP{Assumption III:} $\sup_{t\geq 0} \Expect[(r_t^2+q_t^2)\SigN_t]= M<\infty$. 
\medskip

\begin{proposition}\label{prop:XN}
	Consider the particle system~\eqref{eq:particles} and the mean-field process~\eqref{eq:Xbar} under Assumptions (I)-(II)-(III). Then,  
the error between the empirical variance and the mean-field variance is bounded according to
\begin{align}\label{eq:error}
	\Expect[|\SigN_t-\bar{\Sigma}_t|^2]\leq \text{(const.)} \left[\frac{M}{N} + e^{-2\lambda_0 t}\Expect[|\SigN_0 - \bar{\Sigma}_0|^2]\right]
\end{align}
for $N>\frac{16H^4M}{(\lambda_0-A)^2\lambda_0}$.  
\end{proposition}
This result forms the basis for the convergence of the empirical distribution to the mean-field distribution. However, the analysis is more involved and the subject of ongoing work. 
We conjecture a result of the form
\begin{align} \nonumber
d(\pi^{(N)}_t,&\bar{\pi}_t) \lesssim \frac{M}{N} + e^{\int_0^t G_s \ud s}d(\pi^{(N)}_0,\bar{\pi}_0) \\&+ e^{-2\lambda_0t}(\Expect[|\mN_0-\mbar_0|^2 + |\SigN_0-\Sigmabar_0|^2]) \label{conjecture}
\end{align}
where $d(\mu,\nu) = \sup_{\|\nabla f\|_\infty<1 }\Expect[\mu(f)-\nu(f)|^2]$ is a metric  between two (possibly random) probability measure $\mu$ and $\nu$, and $\mu(f) := \int f \ud \mu$. 
The convergence analysis is carried out in the literature, for the three forms of the EnKF in Table~\ref{tab:Grqs}, under strong assumption that the linear system is stable, i.e. $A<0$~\citep{delmoral2016stability,bishop2018stability,taghvaei2019design} (see Remark~\ref{remark}).  

The error result~\eqref{eq:error} involves two terms. The second term is due to the error in the initial variance. It converges to zero as $t \to \infty$ for any choice of $G_t,r_t,q_t$, as long as~\eqref{eq:Grq-1} holds, due to the stability of the Kalman-Bucy filter. The second term is due to the stochastic terms present in the particle system. It is proportional to $M$ and converges to zero as $N\to\infty$. 
The value of the constant $M$ depends on the choice for $r_t$ and $q_t$. In particular, for the three forms of the EnKF algorithm, we have
\begin{align*}
M =
 	\begin{cases}
	\Sigma_B \sup_t \Expect[\SigN_t] + H^2 \sup_t \Expect[(\SigN_t)^3] ,  &\text{ (P-EnKF)}, \\
	\Sigma_B \sup_t \Expect[\SigN_t],  &\text{  (S-EnKF)}, \\
	0,  &\text{(D-EnKF)}. 
	\end{cases}
\end{align*}
The following conclusions can be drawn: \\(i) P-EnKF admits larger steady-state error compared to S-EnKF, while it was  shown that it is more stable (see~\eqref{eq:W2-EnKF}). \\(ii) The steady-state error for D-EnKF, in approximating the variance, converges to zero, even for finite $N$.  The reason is that the evolution of the empirical variance is deterministic, and identical to Kalman-filter equation, while for P-EnKF and S-EnKF involve  stochastic terms. \\
(iii) Assuming the conjecture~\eqref{conjecture} is true, the steady-state error for  D-EnKF, in approximating the distribution, is proportional to $d(\pi^{(N)}_0,\pibar_0)$ which depends on the initial conditions of the particles. If the initial particles are sampled randomly from $\pibar_0$, then the error decays as $N^{-1}$. However, if the initial particles are not random samples from $\pibar_0$, the error persists to exists even as $N\to \infty$. This is due to the fact that the mean-field system is not stable, as shown in~\eqref{eq:W2-EnKF}. This is in contrast to P-EnKF and S-EnKF where the steady-state error for the distribution is proportional to $\frac{M}{N}$ independent of the initial condition of the particles.   



\begin{remark}\label{remark}
	The result of the Prop.~\ref{prop:XN} also holds in vector-case, under additional assumption that the linear system is stable and fully observable ($H$ is full-rank). For detailed analysis of P-EnKF algorithm, see~\citep{delmoral2016stability,bishop2018stability}, and extension to nonlinear setting~\citep{delmoral2017stability,jana2016stability}. Analysis of S-EnKF and D-EnKF appears in~\citep{Amir-CDC2018} and \citep{AmirACC2018} respectively. 
\end{remark}
\section{Conclusion}\label{sec:conclusion}
The paper presents stability and optimality analysis of EnKF algorithms. The central concept is the construction of a mean-field process that is exact, i.e. its time marginal distribution is equal to the filter posterior distribution. Because the exactness does not specify the joint distribution of the marginal distributions, there are infinity many exact mean-field processes (see Prop.~\ref{prop:Xbart}). Stability and accuracy of three forms of mean-field process, that correspond to three forms of EnKF, are studied (see Table~\ref{tab:Grqs}). It is shown that the deterministic EnKF is not stable, while stochastic forms of the EnKF are stable. The stability is stronger for EnKF with larger stochastic input (see Prop.~\ref{prop:stability} and the proceeding discussion). While stochastic forms of EnKF are stable, they admit larger approximation error because of the presence of stochastic terms. The deterministic EnKF is accurate only if the particles are initialized properly, because it does not correct for the initial error~(see Prop.~\ref{prop:XN} and proceeding discussion).    

The analysis that is presented in this paper raises the question about how to design a  mean-field process in the nonlinear and non-Gaussian setting, such that the mean-field process is stable.

\appendix

\section{Proof of Propisition~\ref{prop:stability}}

(i) 
To obtain a bound for the error in variance, we use the explicit solution to the Riccati equation~\eqref{eq:KF-var} that is available for the scalar case. In particular, let $\phi_t(x)$ denote the semigroup associated with the Riccati equation  such that $\Sigma_t = \phi_t(\Sigma_0)$. The explicit form of the semigroup is given by
	\begin{equation}\label{eq:Riccati-semigroup}
	\phi_t(x) =\frac{(\lambda_0 +A\tanh(\lambda_0 t) )x  + 
		\Sigma_B \tanh(\lambda_0 t)}{\lambda_0 -A\tanh(\lambda_0 t) + H^2 \tanh(\lambda_0 t)x}
	\end{equation} 
	The  bound is obtained by providing a bound for the derivative of the semigroup, with respect to $x$, denoted by $\nabla \phi_t(x)$:
	\begin{align}\nonumber
	\nabla 	\phi_{t}(x)  &= \frac{\lambda_0^2}{\cosh(\lambda_0 t)^2(\lambda_0 -A\tanh(\lambda_0 t) + H^2 \tanh(\lambda_0 t)x)^2} \\&\leq \frac{4\lambda_0^2}{(\lambda_0-A)^2}e^{-2\lambda_0 t}\label{eq:grad-phi}
	\end{align}
	Therefore, $\phi_t$ is globally Lipschitz and 
	\begin{align*}
	|\bar{\Sigma}_t - \Sigma_t|&=|\phi_t(\bar{\Sigma}_0)- \phi_t(\Sigma_0)|\leq c_1e^{-2\lambda_0 t} |\bar{\Sigma}_0-\Sigma_0|
	\end{align*}
	where $c_1 = \frac{4\lambda_0^2}{(\lambda_0-A)^2}$.

To prove the bound for the mean, we subtract the equation for $m_t$ from the equation for $\bar{m}_t$ to obtain
	\begin{align*}
	\ud \bar{m}_t - \ud m_t &= (A-\bar{\Sigma}_tH^2)(\bar{m}_t - m_t)\ud t \\&+ (\bar{\Sigma}_t - \Sigma_t)H\ud I_t
	\end{align*}   
	where $\ud I_t:=\ud Z_t - H\bar{m}_t\ud t$ is the innovation process. Therefore,  the difference $\bar{m}_t - m_t$ solves a linear system for which the solution is given by
		\begin{align*}
	 \bar{m}_t - m_t &= \psi_t (\bar{m}_0 - m_0) + \int_0^t \psi_{t,s} (\bar{\Sigma}_t - \Sigma_t)H\ud I_t
	\end{align*}   
	where $\psi_{t,s} = \exp(\int_s^t (A-\bar{\Sigma}_\tau H^2) \ud \tau)$. Next we obtain  a bound for $\psi_{t,s}$. 
	\begin{align*}
	\psi_{t,s} &= 
	\exp(\int_s^t (A-\Sigma_\infty H^2) \ud \tau + \int_s^t(\Sigma_\infty -\Sigma_\tau)H^2 \ud \tau  )\\
	&\leq  \exp(-(t-s)\lambda + c_1|\bar{\Sigma}_0 - \Sigma_\infty|\int_s^t e^{-2\lambda\tau} \ud \tau )\\
	&\leq c_2 \exp(-(t-s)\lambda )
	\end{align*} 
	where $c_2 = e^{\frac{c_1|\Sigma_0-\Sigma_\infty|}{2\lambda }}$.
	Therefore,
			\begin{align*}
	\Expect[|\bar{m}_t - m_t|^2] &\leq c_2^2 e^{-2\lambda t} (\bar{m}_0 - m_0)^2 \\&+ \int_0^t c_2^2e^{-2\lambda (t-s)}|\bar{\Sigma}_s - \Sigma_s|^2H^2\ud t\\
	&\leq c_2^2 e^{-2\lambda t} (\bar{m}_0 - m_0)^2 + c_3^2e^{-2\lambda t}|\bar{\Sigma}_0 - \Sigma_0|^2
	\end{align*}   
	where $c_3 = \frac{c_2c_1H}{\sqrt{2\lambda }}$. Combining the bounds for the variance and the mean, we obtain the bound for the second term.

	(ii) The result follows from a coupling argument. Introduce a new process $\tilde{X}_t$ with the initial distribution $N(\bar{m}_0,\Sigmabar_0)$ governed by the same equation as the mean-field process:
\begin{align*}
	\ud \tilde{X}_t =  ~&\mathcal T_t(\mbar_t,\Sigmabar_t) + G_t(\tilde{X}_t-\mbar_t)\ud t + r_t \ud \Bbar_t \\\label{eq:Xbar}
	+ &q_t \ud \Wbar_t ,\quad \Xbar_0 \sim  N(\bar{m}_0,\Sigmabar_0)
\end{align*} 
where $\Bbar_t$ and $\Wbar_t$ are the same as in~\eqref{eq:Xbar}. It is straightforward to verify $\tilde{X}_t \sim \tilde{\pi}_t=N(\bar{m}_t,\Sigmabar_t)$. Subtracting the equation for $\tilde{X}_t$ from $\bar{X}_t$ yields 
$
 	\ud (\Xbar_t - \tilde{X}_t) = G_t(\Xbar_t - \tilde{X}_t) \ud t
$
concluding
 \begin{align*}
\Xbar_t - \tilde{X}_t = e^{\int_0^t G_s \ud s} (\Xbar_0 - \tilde{X}_0) 
\end{align*}
The result follows from definition of the Wasserstein distance and the optimal coupling of the initial conditions.\\
(iii) 	The  result follows from triangle inequality:
\begin{equation*}
	\mathcal W_2(\bar{\pi}_t,\pi_t)\leq 	\mathcal W_2(\bar{\pi}_t,\tilde{\pi}_t) + \mathcal W_2(\tilde{\pi}_t,\pi_t) 
\end{equation*} 
and application of part (i), part (ii), and of the formula for Wasserstein distance between two Gaussian distributions. ,
\begin{equation*}
	\mathcal W_2(\tilde{\pi}_t,\pi_t)^2 = (\bar{m}_t-m_t)^2 + (\sqrt{\Sigmabar_t} - \sqrt{\Sigma_t})^2.
\end{equation*} 

\section{Proof of the proposition~\ref{prop:XN}}
The time evolution of the empirical variance $\SigN_t$ is given by:
\begin{equation}\label{eq:SigN}
	\ud \Sigma^{(N)}_t = \text{Ricc}(\Sigma^{(N)}_t) \ud t + \ud \zeta^{(N)}_t 
\end{equation}
where \[ \ud \zeta^{(N)}_t = \frac{2}{N-1}\sum_{i=1}^N (X^i_t-m^{(N)}_t) (r_t\ud B^i_t+q_t \ud W^i_t).\] Note that $\zeta^{(N)}_t$ is a martingale and $\langle \ud \zeta^{(N)}_t\rangle^2= \frac{4(r_t^2+q_t^2)}{N-1}\Sigma^{(N)}_t\ud t$.  
The difference $\Sigma^{(N)}_t - \Sigma_t$ can be expressed, in terms of the semigroup~\eqref{eq:Riccati-semigroup}, according to
\begin{align*}
	\Sigma^{(N)}_t& - \Sigma_t =  \phi_{0}(\Sigma^{(N)}_t) - \phi_t(\Sigma^{(N)}_0) + \phi_t(\Sigma^{(N)}_0)  - \phi_{t}(\Sigma_0) 
\end{align*}
The bound for the second term is straightforward: \[ |\phi_t(\Sigma^{(N)}_0)  - \phi_{t}(\Sigma_0) |\leq c_1e^{-2\lambda_0t} |\SigN_0 - \Sigma_0|.\] 
The first term is 
\begin{align*}
	  \phi_{0}&(\Sigma^{(N)}_t) - \phi_t(\Sigma^{(N)}_0) = \int_0^t \ud_s \phi_{t-s}(\Sigma^{(N)}_s)  
	\\&= \int_0^t \nabla \phi_{t-s}(\Sigma^{(N)}_s)\ud \zeta^{(N)}_s+\int_0^t\frac{1}{2}\nabla^2 \phi_{t-s}(\Sigma^{(N)}_s)\langle \ud \zeta^{(N)}_s\rangle^2
\end{align*}
where $\nabla^2 \phi_t(x)$ denotes the second-order derivative with respect to $x$, and we used $\frac{\ud}{\ud t} \phi_{t}(x) =  \nabla \phi_{t}(x) \Ricc(x)$ and~\eqref{eq:SigN}. Therefore, using the bounds~\eqref{eq:grad-phi} and
\begin{align*}
	|\nabla^2 \phi_t (x)| &
\leq  c_4 e^{-2\lambda_0 t}
\end{align*}
with $c_4=\frac{8\lambda_0^2H^2}{(\lambda_0-A)^3}$, yields 
\begin{align*}
	\Expect[|\phi_{0}&(\Sigma^{(N)}_t) - \phi_t(\Sigma^{(N)}_0)|^2] = \int_0^t \Expect[\nabla \phi_{t-s}(\SigN_s)^2\langle \ud \zeta^{(N)}_s\rangle^2]\\
	&+\left[\int_0^t\Expect[\nabla^2 \phi_{t-s}(\SigN_s) \langle \ud \zeta^{(N)}_s\rangle^2] \right]^2 \\&\leq \frac{4c_1^2}{N-1}\int_0^te^{-4\lambda_0(t-s)}\Expect[(q_s^2+r_s^2) \SigN_s]\ud s \\
	&+ \left[\frac{4c_4}{N-1}\int_0^te^{-2\lambda_0(t-s)}\Expect[(q_s^2+r_s^2) \SigN_s]\ud s  \right]^2\\
	&\leq \frac{2c_1^2M}{\lambda_0N}
\end{align*}
if $N>\frac{4c_4^2M}{c_1^2\lambda_0}=\frac{16H^4M}{(\lambda_0-A)^2\lambda_0}$ and $M := \sup_{t} \Expect[(q_t^2+r_t^2) \SigN_t]$. This concludes the bound:
\begin{align*}
	\Expect[|\SigN_t-\Sigma_t|^2]\leq \frac{4c_1^2M}{\lambda_0N} + 2c_1^2e^{-2\lambda_0 t}\Expect[|\SigN_0 - \Sigma_0|^2]
\end{align*}

\bibliographystyle{plain}
\bibliography{ref,references}

\end{document}